\documentclass[a4paper,runningheads]{llncs}

\usepackage{graphicx}
\usepackage{makeidx}
\usepackage{amsfonts,amssymb
}
\usepackage{enumerate}
\usepackage{amsmath}

\usepackage{url}
\usepackage[colorlinks=true]{hyperref}
\hypersetup{
pdftitle=A non-autonomous stochastic discrete time system with uniform disturbances, 
pdfauthor=Ioannis K. Dassios and Krzysztof J. Szajowski, 
urlcolor=blue,
filecolor=magenta,
citecolor=green, 
linkbordercolor={1 1 1}, 
citebordercolor={1 1 1}, 
urlbordercolor={ 1 1 1}}

\urldef{\mailsa}\path|Krzysztof.szajowski@pwr.edu.pl|
\newcommand{\keywords}[1]{\par\addvspace\baselineskip
\noindent\keywordname\enspace\ignorespaces#1}

\def\rightbox{\protect\vspace*{-2ex}
\begin{flushright}\(\blacksquare\)\end{flushright}}
\newenvironment{ProofIFIP}{{\sc Proof.}\hspace{1mm}}{\rightbox} 

\def\one{{\mathbb I}} 
\def\rank{{\mbox{rank}}}
\def\colspan{{\mbox{colspan}}} 

\mainmatter  

\title{A non-autonomous stochastic discrete time system with uniform disturbances}

\titlerunning{A non-autonomous stochastic system}

\author{Krzysztof Szajowski%
}
\authorrunning{K.Szajowski}

\institute{Institute of Mathematics and Computer Science, Wroc\l{}aw University of Technology,\\
 Wybrze\.ze Wyspia\'{n}skiego 27, 50-370 Wroc\l{}aw, Poland\\
\mailsa\\
\url{http://www.im.pwr.wroc.pl/\~szajow}}

\toctitle{A non-autonomous stochastic system}
\tocauthor{K.Szajowski}
\begin{document}
\title{A non-autonomous stochastic discrete time system with uniform disturbances}
\titlerunning{A non-autonomous stochastic systems}  
%
\author{Ioannis K. Dassios\inst{1}\inst{2} \and  Krzysztof J. Szajowski\inst{3}
}
\authorrunning{I. K. Dassios and K. J. Szajowski} 
%
\tocauthor{Ioannis K. Dassios and Krzysztof J. Szajowski}
\institute{MACSI, Department of Mathematics \& Statistics, University of Limerick, Ireland,
\and 
ERC, Electricity Research Centre, University College Dublin, Ireland,\\
\email{jdasios@math.uoa.gr}\footnote{WWW home page:\\
\texttt{https://ioannisdassios.wordpress.com/research-visits-invited-talks/}},
\and
Wroc\l{}aw University of Technology, Faculty of Pure and Applied Mathematics\\
Wybrze\.ze Wyspia\'nskiego 27, 50-370 Wroc\l{}aw, Poland\\
\email{Krzysztof.Szajowski@pwr.edu.pl}\footnote{WWW home page: http://neyman.im.pwr.edu.pl/\homedir szajow/}
}
\vspace{-1cm}

\maketitle
\vspace{-4ex}

\begin{abstract}
The main objective of this article is to present Bayesian optimal control over a class of non-autonomous linear stochastic discrete time systems with disturbances belonging to a family  of the one parameter uniform distributions. It is proved that the Bayes control for the Pareto priors is the solution of a linear system of algebraic equations. For the case that this linear system is singular, we apply optimization techniques to gain the Bayesian optimal control. These results are extended to generalized linear stochastic systems of difference equations and provide the Bayesian optimal control for the case where the coefficients of these type of systems are non-square matrices. The paper extends the results of the authors developed for system with disturbances belonging to the exponential family.

\end{abstract}
\vspace{-2ex}

\section{\label{intro} Introduction}\vspace{-2ex}
Linear stochastic discrete time systems (or linear matrix stochastic difference equations), are systems in which the variables take their value at instantaneous time points. The horizon of control depends on the problem. The state at instance $n$ depends on random disturbance and the chosen controls. Discrete time systems differ from continuous time ones in that their signals are in the form of sampled data. With the development of the digital computer, the stochastic discrete time system theory plays an important role in the control theory. In real systems, the discrete time system often appears when it is the result of sampling the continuous-time system or when only discrete data are available for use. The investigation aims are, when such system is under consideration, determining the control goals, performance measures and the information available at moments of controls' specification. The small deviations of the parameters can be treated as disturbances. As the random disturbance is admitted the performance measure will be the mean value of the deviation of the states from the required behavior of the system. When all the parameters of the system are known and the distribution of disturbances is well defined then the optimal control can be determined at least for the finite horizon case. The extension of the model to the adaptive one means that the disturbances are not precisely described. Adaptive control is the control method used by a controller which must adapt to a controlled system with parameters which vary, or are initially uncertain (c.f. Black et al.~\cite{blahagary14:adaptive} or Tesfatsion~\cite{Tes82:dual} for the history of the adaptive control). Under some unification the model of adaptive control of the linear system is formulated as a control of the discrete time Markov process (cf. \cite{DunPasSte98:adaptive}).

It is assumed that the disturbance has a fixed probabilistic description which is determined by the assumption. In this paper it is assumed that the distribution function is known to be an accuracy of parameters and the disturbances additionally change the state of the system. It resembles the statistical problem of estimation. It was seminal paper by Wald~\cite{wal39:estimation} where the background of the modern decision theory was established (cf. \cite[Chapt. 7]{wal50:Statistical}). The decision theory approach to the control problems were immediately applied (see books by  Sworder~\cite{Swo66:Adaptive}, Aoki~\cite{Aoki67:Optimization}, Sage and Melsa~\cite{SagMel71:Estimation}). The new class of control systems under uncertainty was  called \emph{adaptive} (cf. \cite{Tes82:dual}, \cite{blahagary14:adaptive}). In these adaptive control problems the important role have Bayesian systems. In this class of control models it is assumed that the preliminary knowledge of  the disturbances is given by \emph{a priori} distributions of their parameters. 
The aim is to construct the controls in a close form. The construction of the Bayes control is also auxiliary for the construction of minimax controls (see Szajowski and Trybu{\l}a~\cite{SzaTry87:minimax}, Porosi\'nski and Szajowski~\cite{PorSza88:mincontrol},  Grzybowski~\cite{Grz91:Minimax}, Gonz{\'a}lez-Trejo et al. ~\cite{GonHerHoy02:Minimax}). It is observed the interested in various models of disturbance structure (cf. Duncan \& Pasik-Duncan~\cite{DunPas13:discrete}) and the disturbance distributions (cf. Walczak~\cite{Wal86short:Bayes,MA827}). Stochastic discrete time systems have many applications which we have described in \cite{DasSza15:Bayes} where the Bayes control of the linear system with quadratic cost function and the disturbances having the distribution belonging to the exponential family with conjugate priors is solved.            

The paper is organized as follows: the description of the stochastic discrete time systems is subject of the section~\ref{subsec1.1} and  some remarks on disturbances are given in the section~\ref{sec2}. In the section~\ref{BCont} we determine the Bayes control for the conjugate prior distribution $\pi$ of the parameter $\overline{\lambda}$ as the solution of a singular linear system and provide optimal Bayesian control. We close the paper by studying the Bayes control of a class of generalized linear stochastic discrete time systems.  
\nopagebreak
\subsection{\label{subsec1.1} Stochastic discrete time systems}\vspace{-2ex} Let $\overline{x}_n\in \mathbb{R}^m$ be the state of the system, $\overline{u}_n\in\mathbb{R}^m$ be the control. Assume that $\overline{\upsilon}_{n}\in V\subset\mathbb{R}^m$, with $\overline{\upsilon}_{n}=(\upsilon_{n}^{1},\upsilon_{n}^{2},\ldots,\upsilon_{n}^{k},0,\ldots,0)^{T}$, is the disturbance at time $n$ and $\alpha_n$, $b_n$, $c_n$$\in\mathbb{R}^{m\times m}$.  Consider a stochastic discrete time system (cf.  Kushner~\cite{Kus71:IntroSC})
\begin{equation}\label{eq1}
\overline{x}_{n+1}=\alpha_n \overline{x}_n + b_n \overline{u}_n + c_n \overline{\upsilon}_n,\quad \forall n=0,1,\ldots,N-1.
\end{equation}
The horizon $N$ of the control, the time up to which the system is controlled, is a random variable, independent of the disturbances $\overline{\upsilon}_0, \overline{\upsilon}_1,\ldots$, and has the following known distribution
\begin{equation}\label{eq2}
P\lbrace N=k\rbrace=p_k,\quad \forall k=0,1,...,M,\quad \sum^{M}_{i=0}p_k=1,\quad p_M\neq 0.
\end{equation}
In the authors paper \cite{DasSza15:Bayes} it was considered the family of the exponentially distributed disturbances. Let us assume here that the disturbances $\upsilon_{n}^{i}$ have the uniform distributions on $[0,\lambda_i]$ with parameter $\lambda_{i}\in\Re^{+}$, $i=1,2,\ldots,k$ and
\[
X_{n}=(\overline{x}_{0},\overline{x}_{1},...,\overline{x}_{n}),\quad U_{n}=(\overline{u}_{0},\overline{u}_{1},...,\overline{u}_{n}), \quad \overline{\lambda}=(\lambda_{1},\lambda_{2},...,\lambda_{k},0,...,0)^{T}.
\]
For convenience $U_{M}$ will be denoted by $U$ and called a control policy. 
\begin{definition}\label{ControlCost} The control cost for a given policy $U$ (the loss function) is the following
\begin{equation}\label{L}
L(U,X_{N})=\sum^{N}_{i=0}(\overline{y}_{i}^{T}s_{i}\overline{y}_{i}+\overline{u}_{i}^{T}k_{i}\overline{u}_{i}),
\end{equation}
where $k_i\in\mathbb{R}^{m\times m}\geq 0_{m,m},$ are symmetric matrices, $s_i\in\mathbb{R}^{2m\times 2m}\geq 0_{2m,2m}$ and $\overline{y}_i=\left(\begin{array}{c}\overline{x}_i\\\cdots\\\overline{\lambda}\end{array}\right)\in\mathbb{R}^{2m},$ $\forall i=0,1,...,M$. With $0_{i,j}$ we will denote the zero matrix $i\times j$.
\end{definition}

Let the prior distribution $\pi$ of the parameter $\overline{\lambda}$ be given. It is considered the Pareto priors (see \cite[Ch. 9.7]{deG70:Optimal}, \cite{Phi85:Pareto}) with parameters $r_i>0$, $\beta_i>2$
\begin{equation}\label{ParetoPrior}
g(\lambda_i|\beta_i,r_i)=\frac{\beta_i r_i^{\beta_i}}{\lambda_i^{\beta_i+1}}\one_{[r_i,\infty)}(\lambda_i).
\end{equation}
Denote $E_{N},E_{\overline{\lambda}}$ the expectations with respect to the distributions of $N$ and random vectors $\overline{\upsilon}_{0},\overline{\upsilon}_{1},...$ (when $\overline{\lambda}$ is the parameter), $E_{\pi}$ and $E$ are the expectations with respect to the distribution $\pi$ and to the joint distribution $\overline{\upsilon}_{n}$ and $\overline{\lambda}$, respectively. 

\begin{definition}\label{RiskFunction} ({\small see \cite{Kus71:IntroSC}, \cite{PorSzaTry85:multi}, \cite{Saw81:Exact}, \cite{Wal86short:Bayes,MA827}}) Let $L(\cdot,\cdot)$ be the loss function given by \eqref{L}. 
\begin{description}
\item[(a)] The risk connected with the control policy $U$, when the parameter $\overline{\lambda}$ is given, is defined as follows
\[
R(\overline{\lambda},U)=E_{N}\left[E_{\overline{\lambda}}[L(U,X_{N})\mid X_0]\right]=E_{N}\left[E_{\overline{\lambda}}[\sum^{N}_{i=0}\overline{y}_{i}^{T}s_{i}\overline{y}_{i}+\overline{u}_{i}^{T}k_{i}\overline{u}_{i}\mid X_0]\right].
\]
\item[(b)] The expected risk $r$, associated with $\pi$ and the control policy $U$, is equal to
\[
r(\pi,U)=E_{\pi}[R(\overline{\lambda},U)]=E_{N}\left[E[\sum^{N}_{i=0}\overline{y}_{i}^{T}s_{i}\overline{y}_{i}+\overline{u}_{i}^{T}k_{i}\overline{u}_{i}\mid X_0]\right].
\]
\item[(c)] The expected risk $r$, associated with $\pi$ and the control policy $U$, is equal to
\[
r_n(\pi,U^n)=E_N[E_{\pi}[R(\overline{\lambda},U)]=E_{N}\left[E[\sum^{N}_{i=n}\overline{y}_{i}^{T}s_{i}\overline{y}_{i}+\overline{u}_{i}^{T}k_{i}\overline{u}_{i}\mid X_0]\right].
\]
\end{description}
\end{definition}

Let the initial state $\overline{x}_{0}$ and the distribution $\pi$ of the parameter $\overline{\lambda}$ be given. 
\begin{definition}\label{def_BayesRisk} 
 A control policy $U^{*}$ is called the Bayes policy when $r(\pi,U^{*})=\inf_{U\in \wp_{\pi}} r(\pi,U)$, 
where $\wp_{\pi}$ is the class of the control policies $U$ for which exists $r(\pi,U)$.
\end{definition}


\subsection{\label{sec2} Filtering}\vspace{-2ex}
Let us assume that the random variables $\overline{\upsilon}_n$ have the density $p(\overline{\upsilon}_{n},\overline{\lambda})$ with respect to a $\sigma$-finite measure $\mu$ on $\mathbb{R}$. The consideration is focused on the special case when each coordinate has the uniform distribution, i.e. the density $p(\overline{\upsilon}_{n},\overline{\lambda})$ has the following representation:
\begin{equation}\label{eq3}
p(\overline{\upsilon}_{n},\overline{\lambda})=\prod^{k}_{i=1}p(\upsilon_{n}^{i},\lambda_{i}),
\end{equation}
where 
$p(\upsilon_{n}^{i},\lambda_{i})=\frac{1}{\lambda_i}\one_{[0,\lambda_i]}(\upsilon_{n}^{i})$, for all $i=1,2,...,k$. $V_{i}^{*}$ is the set of the random variables $\upsilon_n^i$. 
We have: 
\[
E_{\lambda_i}[\upsilon_{n}^{i}]=\frac{\lambda_i}{2}=q_{i}\lambda_{i} \mbox{ and, }
E_{\lambda_i}[(\upsilon_{n}^{i})^2]=\frac{\lambda_i^2}{3}=q_{1,i}\lambda_{i}^{2},
\]
where $q_{i},q_{1,i}$ are constants. Let $\overline{\lambda}$ have the a priori distribution $\pi$ with density
\begin{equation}\label{eq4}
g(\overline{\lambda}\mid \, \overline{\beta},\overline{r})=\prod^{k}_{i=1}g_{i}(\lambda_{i}; \,\beta^{i},r^{i}),
\end{equation}
where $g_{i}(\lambda_{i}| \, \beta^{i},r^{i})$ is given by \eqref{ParetoPrior} where $\overline{\beta}\in S_k^\beta\subset\mathbb{R}^m$, $\overline{r}\in S_k^r\subset\mathbb{R}^m$ with $$\overline{\beta}=(\beta^{1},\beta^{2},\ldots,\beta^{k},0,\ldots,0)^{T},$$ and $\overline{r}=(r^{1},r^{2},\ldots,r^{k},0,\ldots,0)^{T}$. 
When such the \emph{a priori} density is assigned to $\lambda_i$ and then the object of filtering, to determine the Bayes control, is to produce \emph{a posteriori} density for $\lambda_i$ after any new observations of the state of the system. We change the control after obtaining the new data. 
Hence, to determine the Bayes control, a posteriori density for $\overline{\lambda}$ must be obtained after any new observations. This is possible if for $n=0,1,...N-1$ and a given $\overline{x}_0$, we can derive $\overline{\upsilon}_n$ from \eqref{eq1}, i.e. the equations
\[
\overline{\upsilon}_{n}=c_n^{-1}[\overline{x}_{n+1}-\alpha_{n}\overline{x}_{n}-b_{n}\overline{u}_{n}].
\]
If for a value of $n$, the matrix $c_n$ is singular, we will have to compute the Moore-Penrose Pseudoinverse $c^\dagger_n$ and then use the following expression
\[
\overline{\upsilon}_{n}=c_n^{^\dagger}[\overline{x}_{n+1}-\alpha_{n}\overline{x}_{n}-b_{n}\overline{u}_{n}].
\] 
The Moore-Penrose pseudo-inverse can be calculated via the singular value decomposition of $c_n$ (see \cite{Rug93:Linear}). In these cases \emph{a posteriori} density $f(\overline{\lambda}\mid X_{n},U_{n-1})$ of the parameter $\overline{\lambda}$, after having observed $X_{n}$ and chosen $U_{n-1}$, has the same form as \eqref{eq4} i.e.
\[
f(\overline{\lambda}\mid X_{n}, \, U_{n-1})=f(\overline{\lambda}\mid V_{n-1})=g(\overline{\lambda}\mid \, \overline{\beta}_{n},\overline{r}_{n}),
\]
where $V_{n-1}=(\overline{\upsilon}_{0},\overline{\upsilon}_{1},...,\overline{\upsilon}_{n-1})$, $\overline{\beta}_{n}=\overline{\beta}_{n-1}+2\overline{q}$, $\overline{q}=(q_{1},q_{2},...,q_k,0,...,0)^{T}\in Q_i^*\subset\mathbb{R}^m$ and $\overline{r}_{n}=\overline{r}_{n-1}\vee \overline{\upsilon}_{n}$ ($\overline{r}_0=\overline{r}$). Under these denotations we have $E(\lambda_{i}\mid X_{n},U_{n-1})=T^{n,i}r_{n}^{i}=\frac{\beta_n}{\beta_n-1}r_{n}^i$ and $E(\lambda_{i}^{2}\mid X_{n},U_{n-1})=T_{1}^{n,i}(r_{n}^{i})^{2}$. For known $X_{n}$ and $U_{n-1}$, the conditional distribution of $\overline{\upsilon}_{n}$ has the density
\begin{eqnarray*}
h(\overline{\upsilon}_n\mid X_{n},U_{n-1})&=&\prod_{i=1}^{k}h_{i}(\upsilon^i_n\mid X_{n},U_{n-1}),
\end{eqnarray*}
where  
\begin{eqnarray*}\label{densH}
h_{i}(\upsilon^{i}_n\mid X_{n},U_{n-1})&=&\int_0^\infty p(\upsilon^{i}_n,\lambda)g(\lambda|\beta_n^i,r_n^i)d\lambda\\
&=&\frac{\beta_{n}^{i}(r_{n}^{i})^{\beta_n^i}}{\beta_{n+1}^{i}}\frac{1}{(r_{n+1}^{i})^{\beta_{n+1}}}\one_{[0,\infty)}(v),
\end{eqnarray*}
 for $n=0,1,...,M-1$, $i=1,2,...,k$. In addition (see \cite{Wal86short:Bayes,MA827}) by direct calculation we get
\begin{lemma}\label{ExpValues}
The following  equations are fulfilled:
\begin{eqnarray}\label{Elv}
E(v_n^{i}&\mid X_{n},U_{n-1})=&\frac{1}{2}\frac{\beta_n^i }{\beta_{n+1}^{i}}r_{n}^{i}=Q^{n,i}r_n\\
\label{E2v2}
E((v_n^{i})^{2}&\mid X_{n},U_{n-1})=&Q_{1}^{n,i}(r_{n}^{i})^{2}\mbox{where $Q_{1}^{n,i}=\frac{\beta_n^i}{3(\beta_n-2)}$.}\\
\label{E3r}
E(r_{n+1}^{i}&\mid X_{n},U_{n-1})=&Q_2^{n,i}r_{n}^{i} \mbox{ where $Q_2^{n,i}=\frac{(\beta_n^i)^2}{(\beta_n^i)^2-1}$,}\\
\label{E4r2}
E((r_{n+1}^{i})^{2}&\mid X_{n},U_{n-1})=&Q_{3}^{n,i}(r_{n}^{i})^{2}\mbox{where $Q_{3}^{n,i}=\frac{\beta_n^i(\beta_n^i-1)}{(\beta_n+1)(\beta_n+2)}$.}\\
\label{E5x}
E(x_{n+1}^{i}&\mid X_{n},U_{n-1})=&\alpha_nx_n+u_n+\gamma_nQ^{n,i}r_{n}^{i},\\
\label{E6x2}
E((x_{n+1}^{i})^{2}&\mid X_{n},U_{n-1})=&(\alpha_n x_n+u_n)^2+2(\alpha_n x_n+u_n)\gamma_n^i Q^{n,i} r_n^i\\
\nonumber&&\vspace{10em}+\gamma_n^2 Q_1^{n,i}(r_{n}^{i})^{2},\\
\label{E7rx}
E(x_{n+1} r_{n+1}^{i}&\mid X_{n},U_{n-1})=&(\alpha_n x_n+u_n)Q_{2}^{n,i}r_{n}^{i}+\gamma_n^i Q_4^{n,i} r_n^i, 
\end{eqnarray}
where $Q_{4}^{n,i}=\frac{(\beta_n^i)^2}{(\beta_n+1)(\beta_n-2)}$.
\end{lemma}

\section{\label{BCont}The Bayes Control}\vspace{-2ex}
Suppose the initial state $\overline{x}_{0}$ is given, the disturbances have the distribution with the density given by \eqref{eq3} and the prior distribution $\pi$ of the parameter $\overline{\lambda}$ is given by \eqref{eq4}. Let the distribution of the random horizon $N$ be given by \eqref{eq2}. Consider the problem of the Bayes control for the system \eqref{eq1} with the starting point at the moment $n$, when $X_{n}$, $U_{n-1}$ are given. The expected risk is then given by (c.f Defintion~\ref{RiskFunction} (c); see  \cite{PorSzaTry85:multi}, \cite{DasSza15:Bayes})
\begin{equation}\label{rUn}
r_n(\pi,U^n)=E\left[\sum^{M}_{i=n}(\overline{y}_{i}^{T}s_{i}\overline{y}_{i}+\overline{u}_{i}^{T}k_{i}\overline{u}_i)\mid X_n,U_{n-1}\right].
\end{equation}
Let us denote $\varphi_{k}=\sum^{M}_{i=k}p_{i}$. We have
\begin{eqnarray*}
r_{n}
&=&E[\sum^{M}_{i=n}\frac{\varphi_{i}}{\varphi_{n}}(\overline{y}_{i}^{T}s_{i}\overline{y}_{i}+\overline{u}_{i}^{T}k_{i}\overline{u}_{i})\mid X_{n},U_{n-1}].
\end{eqnarray*}

For the above truncated problem we provide the following definitions:
\begin{definition}\label{Def3.2}
The Bayes risk is defined as
\begin{equation}\label{W1}
W_{n}=\inf_{U^n} r_n(\pi,U^n),
\end{equation}
where $r(\pi,U^n)$ is the expected risk defined in the defintion~\ref{RiskFunction} (c) and the formulae~\eqref{rUn}.
\end{definition}
\begin{definition}\label{Def3.3}
If there exists $U^{n^{*}}=(\overline{u}_{n}^{*},\overline{u}_{n+1}^{*},..,\overline{u}_{N}^{*})$ such that 
$W_{n}=r(\pi,U^{n^*})$, then $U^{n^{*}}$ will be called the Bayes policy and $\overline{u}_{i}^{*}$, $i=n,n+1,...,N$ the Bayes controls for truncated control problem.
\end{definition}
Obviously,  $r(\pi,U^{0})=r(\pi,U), \quad W_{0}=r(\pi,U^{*})$. For the solution of the Bayes control problem we derive the Bayes controls $\overline{u}_{n}^{*}$ for $n=N,N-1,...,1,0$ recursively. Then $U^{0^{*}}$ is the solution of the problem. From the Bellman's dynamic programming optimality principle we obtain the following Lemma, see \cite{PorSzaTry85:multi}.

\begin{lemma}\label{Lemma3.1} Assume the stochastic discrete time system \eqref{eq1}. Then the Bayes risk $W_{n}$ has the form
\begin{equation}\label{W2}
W_{n}=\overline{x}_{n}^{T}A_{n}\overline{x}_{n}+2\overline{r}_{n}^{T}B_{n}\overline{x}_{n}+2\overline{r}_{n}^{T}C_{n}\overline{r}_{n},
\end{equation}
where $A_{n},B_{n},C_{n}\in\mathbb{R}^{m\times m}$, $\overline{D}_{n}\in\mathbb{R}^{m}$ with $A_{n}=f_1(s_n)$, $B_{n}=f_2(Q^n,Q^n_2,s_n)$, $C_{n}=f_3(Q^n,Q^n_1,Q^n_3,Q^n_4,s_n)$. The functions $f_j$, $j=1,2,3$ are strictly monotonic, differentiable. The constants $Q^{n,i}$, $Q_{j}^{n,i}$, $j=1,2,3,4$, $n=0,1,...,N$ are given by \eqref{Elv}, \eqref{E2v2} and $s_n$ is defined in \eqref{L}.
\end{lemma}

\subsection{\label{SDTSquadratic}Bayesian optimal control for stochastic discrete time systems}\vspace{-2ex}
We can now prove the following theorem
\begin{theorem}\label{Th3.1} Assume the stochastic discrete time system \eqref{eq1}. Then, the Bayes control $\overline{u}_{n}^{*}$ is given by the solution of the linear system
\begin{equation}\label{u*}
K_n\overline{u}_{n}^{*}=L_n,
\end{equation}
where
\begin{equation}\label{kn}
K_n=k_n+b_n^TA_{n+1}b_n
\end{equation}
and
\begin{equation}\label{ln}
L_n=-b_n^T[A_{n+1}\alpha_n \overline{x}_n+(A_{n+1}c_nQ^n+B_{n+1}Q^n_2)\overline{r}_n].
\end{equation}
The matrices $k_n$,  $A_n$, $Q^n$, are defined in \eqref{L},  \eqref{W2}, the lemma~ \ref{ExpValues}, respectively and $\overline{e}=\sum_{j=0}^{n-1}Q^j\overline{r}_j$.
\end{theorem}
\begin{ProofIFIP} From \eqref{W1}, the Bayes risk is given by
$W_{n}=\inf_{U^n} r(\pi,U^n)
$. It is, equivalently,
\[
W_{n}=\min_{U^n} E\left[\sum^{M}_{i=n}(\overline{y}_{i}^{T}s_{i}\overline{y}_{i}+\overline{u}_{i}^{T}k_{i}\overline{u}_i)\mid X_n,U_{n-1}\right].
\]
We have 
\begin{eqnarray*}
W_{n}&=&\min_{\overline{u}_n} \big\{\overline{u}_n^{T}k_n\overline{u}_n+E\left[\overline{y}_n^{T}s_n\overline{y}_n\mid X_n,U_{n-1}\right]\\
&&+\min_{U^{n+1}}E\left[E\left[\sum^{k}_{i=n+1}(\overline{y}_{i}^{T}s_{i}\overline{y}_{i}+\overline{u}_{i}^{T}k_{i}\overline{u}_i)\right]\mid X_n,U_{n-1}\right]\big\}.
\end{eqnarray*}
It means $W_{n}=\min_{\overline{u}_n} \left\{\overline{u}_n^{T}k_n\overline{u}_n+E\left[\overline{y}_n^{T}s_n\overline{y}_n\mid X_n,U_{n-1}\right]+ E\left[W_{n+1}\mid X_n,U_{n-1}\right]\right\}$.
Hence, the Bayes control $\overline{u}_{n}^{*}$ satisfies the equation ($\nabla$ is the gradient):
\[
\nabla_{\overline{u}_{n}}\left\{\overline{u}_n^{T}k_n\overline{u}_n+E\left[\overline{y}_n^{T}s_n\overline{y}_n\mid X_n,U_{n-1}\right]+
E\left[W_{n+1}\mid X_n,U_{n-1}\right]\right\}_{\overline{u}_{n}=\overline{u}_{n}^*}=0_{m,1}.
\]
By using \eqref{W2} we get
\begin{eqnarray*}
k_n\overline{u}_{n}^{*}&+&b_n^TA_{n+1}(\alpha_n \overline{x}_n + b_n \overline{u}_n\\
 &+& c_n E(\overline{\upsilon}_{n}\mid X_{n},U_{n-1}))+E\left[\left\{b_n^TB_{n+1}\overline{r}_{n+1}\right\}_{\overline{u}_n=\overline{u}_{n}^*}\mid X_n,U_{n-1}\right]=0_{m,1}.
\end{eqnarray*}
By the properties of conjugate priors for the uniform distribution  (see the lemma~\ref{ExpValues} we have
\[
k_n\overline{u}_{n}^{*}+b_n^TA_{n+1}(\alpha_n \overline{x}_n + b_n \overline{u}_n + c_n Q^{n}\overline{r}_n)+E\left[\left\{b_n^TB_{n+1}\overline{r}_{n+1}\right\}_{\overline{u}_n=\overline{u}_{n}^*}\mid X_n,U_{n-1}\right]=0_{m,1},
\]
and at the end
$(k_n+b_n^TA_{n+1}b_n)\overline{u}_{n}^{*}=-b_n^T[A_{n+1}\alpha_n \overline{x}_n+(A_{n+1}c_nQ^n+B_{n+1}Q_2^n)\overline{r}_n]$. The proof is completed.
\end{ProofIFIP}
Similarly like for the system with the disturbances belonging to the exponential family (see \cite{DasSza15:Bayes}) we get
\begin{theorem}\label{Th3.2} Consider the system \eqref{eq1} and the matrices $K_n$, $L_n$ as defined in \eqref{kn}, \eqref{ln} respectively. Then
\begin{enumerate}[(a)]
\item $\forall n$ such that $K_n$ is full rank, the Bayes control $\overline u_{n}^{*}$, is given by
\begin{equation}\label{s1}
\overline{{u}}_{n}^{*}=K^{-1}_nL_n.
\end{equation}
\item $\forall n$ such that $K_n$ is rank deficient, the Bayesian optimal control $\hat u_{n}^*$ is given by
\begin{equation}\label{s2}
   \hat u_{n}^*=(K_n^TK_n+E^TE)^{-1}K_n^TL_n.
\end{equation}
Where $E$ is a matrix such that $K_n^TK_n+E^TE$ is invertible and $\left\|E\right\|_2=\theta$, $0<\theta<<1$. Where $\left\|\cdot\right\|_2$ is the Euclidean norm.
\end{enumerate}
\end{theorem}

\subsection{Bayesian optimal control for generalized stochastic discrete time systems}\vspace{-2ex}
In this subsection we will expand the results of the section~\ref{SDTSquadratic} by studying Bayesian optimal control for a class of linear stochastic discrete time systems with non-square coefficients. We consider the following non-autonomous linear stochastic discrete time system
\begin{equation}\label{eq2a}
I_{r,m}\overline{x}_{n+1}=\alpha_n \overline{x}_n + b_n \overline{u}_n + c_n \overline{\upsilon}_n,\quad \forall n=0,1,...,N-1.
\end{equation}
Where $\overline{x}_n\in \mathbb{R}^m$ is the state of the system, $\overline{u}_n\in\mathbb{R}^m$ is the control, $\overline{\upsilon}_{n}\in V\subset\mathbb{R}^m$, with $\overline{\upsilon}_{n}=(\upsilon_{n}^{1},\upsilon_{n}^{2},...,\upsilon_{n}^{k},0,...,0)^{T}$, is the disturbance at time $n$ and $\alpha_n$, $b_n$, $c_n$$\in\mathbb{R}^{r\times m}$. The horizon $N$ of the control is fixed and independent of the disturbances $\overline{\upsilon}_n$, $n\geq 0$. If $r=m$, then $I_{r,m}=I_m$. If $r>m$, then $I_{r,m}=\left[\begin{array}{c}I_m\\0_{r-m,m}\end{array}\right]$ and if $r<m$, then $I_{r,m}=\left[\begin{array}{cc}I_r&0_{r,m-r}\end{array}\right]$ with $I_m, I_r$ identity matrices.
\begin{definition}\label{Def3.4} We will refer to system \eqref{eq2a} as a generalized stochastic linear discrete time system.
\end{definition} 
In the above definition we use the term ''\emph{generalized}'' because the coefficients in the system \eqref{eq2a} can be either square or non-square matrices.
\begin{theorem}\label{Th3.3} Consider the system \eqref{eq2a} for $r\neq m$ and assume the matrices $K_n$, $L_n$ as defined in \eqref{kn}, \eqref{ln} respectively. Then, $\forall n$ such that
\begin{enumerate}[(a)]
\item $m< r$, $\rank (K_n)=m$ and $L_n\in \colspan{ K_n}$, the Bayes control $\overline{{u}}_{n}^{*}$, is given by
\begin{equation}\label{i}
\overline{{u}}_n^*=K^{-1}_nL_n.
\end{equation}
\item $m< r$, $\rank(K_n)=m$ and $L_n\notin \colspan(K_n)$, a Bayesian optimal control is given by
\begin{equation}\label{ii}
    \hat{u}_{n}^*=(K_n^TK_n)^{-1}K_n^TL_n.
\end{equation}
\item $L_n\notin \colspan{K_n}$ and $K_n$ is rank deficient, a Bayesian optimal control is given by
\begin{equation}\label{iii}
   \hat{u}_{n}^*=(K_n^TK_n+E^TE)^{-1}K_n^TL_n.
\end{equation}
Where $E$ is a matrix such that $K_n^TK+E^TE$ is invertible and $\left\|E\right\|_2=\theta$, $0<\theta<<1$.
\item $m> r$, $K_n$ is full rank, a Bayesian optimal control is given by 
\begin{equation}\label{iv}
\hat{u}_{n}^*=K_n^T(K_nK_n^T)^{-1}L_n.
\end{equation}
\item $L_n\in \colspan{K_n}$ and $K_n$ is rank deficient, a Bayesian optimal control is given by \eqref{iii}.
\end{enumerate}
\end{theorem}
The proof is based on ideas similar to those used in prove \cite[Th. 3]{DasSza15:Bayes} and is omitted here. 

\section{Conclusions}\vspace{-2ex}
In this article we focused on developing the Bayesian optimal control for a class of non-autonomous linear stochastic discrete time systems of type \eqref{eq1}. Firstly, we proved that the Bayes control of these type of systems is the solution of a linear system of algebraic equations which can also be singular. For this case we used optimization techniques to derive the optimal Bayes control for \eqref{eq1}. In addition, we used these methods to obtain the Bayesian optimal control of the non-autonomous linear stochastic discrete time system of type \eqref{eq2}, where the coefficients of this system are non-square matrices.

The further extension of this paper is to study to Bayes control problem of stochastic fractional discrete time systems. The fractional nabla operator is a very interesting tool when applied to systems of difference equations and has many applications especially in macroeconomics, since it succeeds to provide information from a specific year in the past until the current year. For all these there is some research in progress.
\vspace{-3ex}
\subsection*{Acknowledgments}
I. Dassios is supported by Science Foundation Ireland (award 09/SRC/E1780). 

\end{document}